\documentclass[11pt]{article}

\usepackage{amsfonts, amsmath, amssymb, amscd, amsthm}

\newtheorem{thm}{Theorem}[section]
\newtheorem{cor}[thm]{Corollary}
\newtheorem{lem}[thm]{Lemma}
\newtheorem{prop}[thm]{Proposition}

\theoremstyle{definition}
\newtheorem{defn}[thm]{Definition}
\theoremstyle{remark}

\newcommand{\prhh }{properly relatively hyperbolic }
\renewcommand{\phi }{{\rm\bf Lab\, }}

\newcommand{\G }{\Gamma (G, X\cup \mathcal H)}

\newcommand{\Hl }{\{ H_\lambda \} _{ \lambda \in \Lambda  }}
\newcommand{\Vl }{\{ V_\lambda \} _{ \lambda \in \Lambda  }}

\newcommand{\N}{\mathbb{N}}
\newcommand{\Z}{\mathbb{Z}}
\newcommand{\Ghh}{$\Gamma (G, X\cup \mathcal{H}')$}

\begin{document}

\title{The SQ-universality and residual properties of relatively hyperbolic groups}
\author{G. Arzhantseva  \and A. Minasyan \thanks{The work of the first two authors was supported by the Swiss
National Science Foundation Grant $\sharp$~PP002-68627.}
\and D. Osin \thanks{The work of the third author has been
supported by the Russian Foundation for Basic Research Grant $\sharp$~03-01-06555.}}
\date{}%

\maketitle

\begin{abstract}
In this paper we study residual properties of relatively hyperbolic
groups. In particular, we show that if a group $G$ is
non-elementary and hyperbolic relative to a collection of proper
subgroups, then $G$ is SQ-universal.
\end{abstract}

\section{Introduction}

The notion of a group hyperbolic relative to a collection of
subgroups was originally suggested by Gromov \cite{Gro} and since
then it has been elaborated from different points of view
\cite{Bow,F,DSO,RHG}. The class of relatively hyperbolic groups
 includes many examples. For instance, if $M$ is a complete
finite-volume manifold of pinched negative sectional curvature,
then $\pi _1(M)$ is hyperbolic with respect to the cusp subgroups
\cite{Bow,F}. More generally, if $G$ acts isometrically and
properly   discontinuously on a proper hyperbolic metric space $X$
so that the induced action of $G$ on $\partial X$ is geometrically
finite, then $G$ is hyperbolic relative to  the collection of
maximal parabolic subgroups \cite{Bow}. Groups acting on $CAT(0)$
spaces with isolated flats are hyperbolic relative to the
collection of flat stabilizers \cite{KH}. Algebraic examples of
relatively hyperbolic groups include free products and their small
cancellation quotients \cite{RHG}, fully residually free groups
 (or Sela's limit groups) \cite{Dah}, and, more generally,
groups acting freely  on $\mathbb R^n$-trees \cite{Gui}.

The main  goal of this paper is to study residual properties of
relatively hyperbolic groups. Recall that a group $G$ is called
{\it SQ-universal} if  every countable group can be
embedded into a quotient of $G$ \cite{Sch}. It is straightforward
to see that any SQ-universal group contains an infinitely
generated free subgroup. Furthermore, since the set of all
finitely generated groups is uncountable and  every single
quotient of $G$ contains (at most) countably many finitely
generated subgroups,  every SQ-universal group has uncountably
many non-isomorphic quotients. Thus the property of being
SQ-universal may, in a very rough sense, be considered as an
indication of ``largeness'' of a group.

The first non-trivial example of an SQ-universal group was
provided by Higman, Neumann and Neumann \cite{HNN}, who proved
that the free group of rank $2$ is SQ-universal. Presently many
other classes of groups are known to be SQ-universal:
various HNN-extensions and amalgamated products
\cite{FT,Los,Sas}, groups of deficiency $2$ \cite{BP}, most $C(3)\,
\& \, T(6)$-groups \cite{How}, etc. The SQ-universality of
non-elementary hyperbolic groups was proved by Olshanskii
\cite{Ols} and, independently, by Delzant \cite{Delzant}. On the other hand, for relatively hyperbolic groups,
there are some partial results. Namely, in \cite{Fine} Fine proved
the SQ-universality of certain Kleinian groups. The case of
fundamental groups of hyperbolic 3-manifolds was studied by
Ratcliffe in \cite{Rat}.

In this paper we prove the SQ-universality of relatively
hyperbolic groups in the most general  settings. Let a group $G$
be hyperbolic relative to a collection of subgroups $\Hl $ (called
{\it peripheral subgroups}).  We say that $G$ is {\it properly
hyperbolic relative to $\Hl $} (or $G$ is a {\it PRH group} for
brevity), if $H_\lambda \ne G$ for all $\lambda \in \Lambda $.
Recall that a group is {\it elementary}, if it
contains a cyclic subgroup of finite index. We observe that every
non-elementary PRH group has a unique maximal finite normal
subgroup denoted by $E_G(G)$ (see Lemmas \ref{prop} and
\ref{Ashot's Lemma} below).

\begin{thm}\label{SQ} Suppose that a group $G$ is
non-elementary and \prhh with respect to a collection of
subgroups $\Hl $. Then for  each finitely generated group $R$,
there exists  a quotient group $Q$ of $G$ and an embedding
$R\hookrightarrow Q$ such that:

\begin{enumerate}
\item $Q$ is properly relatively hyperbolic with respect to the
collection $\{ \psi(H_\lambda) \} _{\lambda \in \Lambda } \cup
\{R\}$ where $\psi\colon G \to Q$ denotes the natural epimorphism;

\item  For each $\lambda \in \Lambda$, we have $H_\lambda \cap
\ker(\psi)= H_\lambda \cap E_G(G)$, that is, $\psi(H_\lambda)$ is
naturally isomorphic to $H_\lambda/(H_\lambda \cap E_G(G))$.

\end{enumerate}
\end{thm}

In general, we can not require the epimorphism $\psi$ to be
injective on every $H_\lambda$. Indeed, it is easy to show that a
finite normal subgroup of a relatively hyperbolic group must be
contained in each infinite peripheral subgroup (see Lemma
\ref{lem:qbf}). Thus the image of $E_G(G)$ in $Q$ will have to be
inside $R$ whenever $R$ is infinite. If, in addition, the group
$R$ is torsion-free, the latter inclusion implies $E_G(G) \le
\ker(\psi)$. This would be the case if one took $G=F_2 \times
\Z/(2\Z)$ and $R=\Z$, where $F_2 $ denotes the free group of rank
$2$ and $G$ is properly hyperbolic relative to its subgroup
$\Z/(2\Z)=E_G(G)$.

Since any countable group is embeddable into a finitely generated
group, we obtain the following.

\begin{cor}\label{SQ-cor}
Any non-elementary PRH group is SQ-universal.
\end{cor}

Let us mention  a particular case of Corollary \ref{SQ-cor}. In
\cite{FT} the authors asked whether every finitely generated group with infinite number
of ends is SQ-universal. The celebrated Stallings theorem
\cite{Stall} states that a finitely generated group has infinite
number of ends if and only if it splits as a nontrivial
HNN-extension or amalgamated product over a finite subgroup. The
case of amalgamated products was considered by Lossov who provided
the positive answer in \cite{Los}. Corollary \ref{SQ-cor} allows
us to answer the question in the general case. Indeed, every group
with infinite number of ends is non-elementary and properly
relatively hyperbolic, since the action of such a group on the
corresponding Bass-Serre tree satisfies  Bowditch's definition of
relative hyperbolicity \cite{Bow}.

\begin{cor}
 A finitely generated group with infinite number of ends is SQ-universal.
\end{cor}

The methods used in the proof of Theorem \ref{SQ} can also be applied
to obtain other results:

\begin{thm}\label{comquot}
Any two  finitely generated non-elementary PRH groups $G_1, G_2$
have a common non-elementary PRH quotient $Q$. Moreover, $Q$ can
be obtained from the free product $G_1\ast G_2$ by adding finitely
many relations.
\end{thm}

In \cite{Ols00} Olshanskii proved that any non-elementary
hyperbolic group has a non-trivial finitely presented quotient
without proper subgroups of finite index. This result was used by
Lubotzky and Bass \cite{BL} to construct representation rigid
linear groups of non-arithmetic type thus solving in negative the
Platonov Conjecture. Theorem \ref{comquot} yields a generalization
of Olshanskii's result.

\begin{defn}
Given a class of groups $\mathcal G$, we say that a group $R$ is
{\it residually incompatible with} $\mathcal G$ if for any group
$A\in \mathcal G$, any homomorphism $R\to A$ has a trivial image.
\end{defn}

If $G$ and $R$ are finitely presented groups, $G$ is properly
relatively hyperbolic, and $R$ is residually incompatible with a
class of groups $\mathcal G$, we can apply Theorem \ref{comquot}
to $G_1=G$ and $G_2=R\ast R$.  Obviously, the obtained common
quotient of $G_1$ and $G_2$ is finitely presented and residually
incompatible with $\mathcal G$.

\begin{cor}\label{ResInc}
Let $\mathcal G$ be a class of groups. Suppose that there exists a
finitely presented group $R$ that is residually incompatible with
$\mathcal G$. Then  every finitely presented non-elementary PRH group has a
non-trivial finitely presented quotient group that is residually incompatible
with $\mathcal G$.
\end{cor}

Recall that there are finitely presented groups having no non-trivial
recursively presented quotients with decidable word problem
\cite{Mil}. Applying the previous corollary to the class $\mathcal
G$ of all recursively presented groups with decidable word
problem, we obtain the following  result.

\begin{cor}
Every non-elementary finitely presented PRH group has an infinite finitely presented quotient
group $Q$ such that the word problem is undecidable in each
non-trivial quotient of $Q$.
\end{cor}

In particular, $Q$ has no proper subgroups of finite index. The reader can
easily check that Corollary \ref{ResInc} can also be applied to the
classes of all torsion (torsion-free, Noetherian, Artinian,
amenable, etc.) groups.

\section{Relatively hyperbolic groups}

We recall the definition of relatively hyperbolic groups suggested
in \cite{RHG} (for equivalent definitions in the case of finitely
generated groups see \cite{Bow,DSO,F}). Let $G$ be a group, $\Hl $ a
fixed collection of subgroups of $G$ (called {\it peripheral
subgroups}), $X$ a subset of $G$. We say that $X$ is a {\it relative
generating set of $G$} with respect to $\Hl $ if $G$ is generated by
$X$ together with the union of all $H_\lambda $ (for convenience,
we always assume that $X=X^{-1}$). In this situation the group $G$ can be
considered as a quotient of the free product
\begin{equation}
F=\left( \ast _{\lambda\in \Lambda } H_\lambda  \right) \ast F(X),
\label{F}
\end{equation}
where $F(X)$ is the free group with the basis $X$.
Suppose that $\mathcal R $ is a subset of $F$ such that the kernel
of the natural  epimorphism $F\to G$ is a normal closure of $\mathcal R $ in the group $F$, then
we say that $G$ has {\it
relative presentation}  \begin{equation}\label{G} \langle X,\; \{
H_\lambda\}_{\lambda\in \Lambda} \mid R=1,\, R\in\mathcal R
\rangle .
\end{equation}
If sets $X$ and $\mathcal R$ are finite, the
 presentation (\ref{G}) is said to be {\it relatively finite}.

\begin{defn} We set
\begin{equation}\label{H}
\mathcal H=\bigsqcup\limits_{\lambda\in \Lambda} (H_\lambda
\setminus \{ 1\} ) .
\end{equation}
A group $G$ is {\it  relatively hyperbolic with respect to a
collection of subgroups} $\Hl $, if $G$ admits a relatively finite
presentation (\ref{G}) with respect to $\Hl $ satisfying a {\it linear
relative isoperimetric inequality}. That is, there exists $C>0$
satisfying the following condition. For every word $w$ in the
alphabet $X\cup \mathcal H$ representing the identity in the group
$G$, there exists an expression
\begin{equation}
w=_F\prod\limits_{i=1}^k f_i^{-1}R_i^{\pm 1}f_i \label{prod}
\end{equation}
with the equality in the group $F$, where $R_i\in \mathcal R$,
$f_i\in F $, for $i=1, \ldots , k$, and  $k\le C\| w\| $, where
$\| w\| $ is the length of the word $w$.  This definition is
independent of the choice of the (finite) generating set $X$ and
the (finite) set $\mathcal R$ in (\ref{G}).
\end{defn}

For a  combinatorial path $p$ in the Cayley graph $\G$ of
$G$ with respect to $X\cup \mathcal H$, $p_-$, $p_+$, $l(p)$, and
$\phi (p)$ will denote the initial point, the ending point, the length
(that is, the number of edges) and the label of $p$ respectively.
Further, if $\Omega$ is a subset of $G$ and $g \in \langle \Omega \rangle \le G$, then
$|g|_\Omega$ will be used to denote the length of a shortest word in $\Omega^{\pm 1}$ representing
$g$.

Let us recall some terminology introduced in \cite{RHG}. Suppose $q$ is a
path in $\G $.
\begin{defn}
A subpath $p$ of $q$ is called an {\it $H_\lambda $-component}
for some $\lambda \in \Lambda $ (or simply a {\it component}) of
$q$, if the label of $p$ is a word in the alphabet
$H_\lambda\setminus \{ 1\} $ and $p$ is not contained in a bigger
subpath of $q$ with this property.

Two components $p_1, p_2$ of a path $q$ in $\G $ are called {\it
connected} if they are $H_\lambda $-components for the same
$\lambda \in \Lambda $ and there exists a path $c$ in $\G $
connecting a vertex of $p_1$ to a vertex of $p_2$ such that ${\phi
(c)}$ entirely consists of letters from $ H_\lambda $. In
algebraic terms this means that all vertices of $p_1$ and $p_2$
belong to the same coset $gH_\lambda $ for a certain $g\in G$.
 We can always assume $c$ to have length at most $1$, as
every nontrivial element of $H_\lambda $ is included in the set of
generators.  An $H_\lambda $-component $p$ of a path $q$ is
called {\it isolated } if no distinct $H_\lambda $-component of
$q$ is connected to $p$.  A path $q$ is said to be {\it without
backtracking} if all its components are isolated.
\end{defn}

The next lemma is a simplification of Lemma 2.27 from \cite{RHG}.

\begin{lem}\label{Omega}
Suppose that a group $G$ is hyperbolic relative to a collection of
subgroups $\Hl $. Then there exists a finite subset $\Omega
\subseteq G$ and a constant $K\ge 0$ such that the following
condition holds. Let $q$ be a cycle in $\G $, $p_1, \ldots , p_k$
a set of isolated $H_\lambda $-components of $q$ for some
$\lambda \in \Lambda $,  $g_1, \ldots , g_k$ elements of $G$
represented by labels $\phi(p_1), \ldots , \phi(p_k)$
respectively. Then $g_1, \ldots , g_k$ belong to the subgroup
$\langle \Omega\rangle \le G$  and the word lengths of $g_i$'s
with respect to $\Omega $ satisfy the inequality
$$ \sum\limits_{i=1}^k |g_i|_\Omega \le Kl(q).$$
\end{lem}

\section{Suitable subgroups of relatively hyperbolic groups}

 Throughout this section let $G$ be a group which is properly
hyperbolic relative to a collection of subgroups $\Hl $, $X$ a
finite relative generating set of $G$, and $\G $ the Cayley graph
of $G$ with respect to the generating set $X\cup \mathcal H$,
where $\mathcal H$ is given by (\ref{H}). Recall that an element
$g\in G$ is called {\it hyperbolic} if it is not conjugate to an
element of some $H_\lambda $, $\lambda\in \Lambda $. The following
description of elementary subgroups of $G$ was obtained in
\cite{ESBG}.

\begin{lem}\label{Eg}
Let $g$ be a hyperbolic element of infinite order of $G$. Then the
following conditions hold.
\begin{enumerate}
\item The element $g$ is contained in a unique maximal elementary
subgroup $E_G(g)$ of $G$, where
\begin{equation} \label{eq:elem} E_G(g)=\{ f\in G\; :\;
f^{-1}g^nf=g^{\pm n}\; {\rm for \; some\; } n\in \mathbb N\}.
\end{equation}

\item The group $G$ is hyperbolic relative to the collection
$\Hl\cup \{ E_G(g)\} $.
\end{enumerate}
\end{lem}

Given a subgroup $S\le G$, we denote by $S^0$ the set of all
hyperbolic elements of $S$ of infinite order. Recall that two
elements $f,g\in G^0$ are said to be {\it commensurable} (in G) if
$f^k$ is conjugated to $g^l$ in $G$ for some non-zero  integers
$k$ and $l$.

\begin{defn} A subgroup $S\le G$ is called
{\it suitable}, if there exist at least two non-commensurable
elements $f_1, f_2\in S^0$, such that $E_G(f_1)\cap E_G(f_2)=\{1\}$.
\end{defn}

If $S^0\ne \emptyset $, we define
 $$E_G(S)=\bigcap\limits_{g\in S^0} E_G(g).$$

\begin{lem} \label{Ashot's Lemma}
If $S \le G$ is a non-elementary subgroup and $S^0 \neq \emptyset$,
then $E_G(S)$ is the maximal finite subgroup of $G$ normalized by
$S$.
\end{lem}

\begin{proof}
Indeed, if a finite subgroup $M \le G$ is normalized by $S$, then
$|S:C_S(M)|<\infty$ where $C_S(M)=\{g \in S\; : \; g^{-1}xg=x,
~\forall \; x \in M\}$. Formula \eqref{eq:elem}
implies that $M\le E_G(g)$ for every $g\in S^0$, hence $M \le
E_G(S)$.

On the other hand, if $S$ is non-elementary and $S^0\ne \emptyset
$, there exist $h \in S^0$ and $a \in S^0 \setminus E_G(h)$. Then
$a^{-1}ha \in S^0$ and the intersection $E_G(a^{-1}ha) \cap
E_G(h)$ is finite. Indeed if $E_G(a^{-1}ha) \cap E_G(h)$ were
infinite, we would have  $(a^{-1}ha)^n=h^k$ for some $n,k\in \mathbb
Z\setminus \{ 0\}$, which would contradict to $a\notin E_G(h)$.  Hence
$E_G(S)\le E_G(a^{-1}ha) \cap E_G(h)$ is finite. Obviously,
$E_G(S)$ is normalized by $S$ in $G$.
\end{proof}

The main result of this section is the following

\begin{prop}\label{suit}
Suppose that a group $G$ is hyperbolic relative to a collection
$\Hl $ and $S$ is a subgroup of $G$. Then the following conditions
are equivalent.  \begin{enumerate} \item[(1)] $S$ is suitable;
\item[(2)] $S^0\ne \emptyset $ and $E_G(S)=\{1\}$.
\end{enumerate}
\end{prop}

Our proof of Proposition \ref{suit} will make use of several
auxiliary statements below.

\begin{lem}[Lemma 4.4, \cite{ESBG}]\label{ah}
For any $\lambda \in \Lambda $ and any element $a\in G\setminus
H_\lambda $, there  exists a finite subset $\mathcal F_\lambda
=\mathcal F_\lambda (a) \subseteq H_\lambda $ such that if $h\in
H_\lambda \setminus \mathcal F_\lambda $, then $ah$ is a
hyperbolic element of infinite order.
\end{lem}

 It can be seen from Lemma \ref{Eg} that 
 every hyperbolic element $g \in G$ of infinite
order is contained inside the elementary subgroup
$$E^+_G(g)=\{f \in G\; :\; f^{-1}g^nf=g^{n}\; {\rm for \; some\; } n\in \N\} \le E_G(g),$$ and
$|E_G(g):E^+_G(g)| \le 2$.

\begin{lem} \label{lem:E^+_G} Suppose $g_1,g_2 \in G^0$ are non-commensurable and  $A=\langle g_1, g_2 \rangle \le G$.
Then  there exists an element $h \in A^0$  such
that:
\begin{enumerate}
\item $h$ is not commensurable with $g_1$ and $g_2$;

\item $E_G(h)=E^+_G(h) \le \langle h, E_G(g_1) \cap E_G(g_2)
\rangle$. If, in addition, $E_G(g_j)=E^+_G(g_j)$, $j=1,2$, then
$E_G(h)=E^+_G(h) =\langle h \rangle \times (E_G(g_1) \cap
E_G(g_2))$.
\end{enumerate}
\end{lem}

\begin{proof} By Lemma \ref{Eg}, $G$ is hyperbolic relative to the collection of peripheral subgroups
$\mathfrak{C}_1=\Hl \cup \{E_G(g_1)\}\cup \{E_G(g_2)\}$. The
center $Z(E^+_G(g_j))$ has finite index in $E^+_G(g_j)$, hence
(possibly, after replacing $g_j$ with a power of itself) we can
assume that $g_j \in Z(E^+_G(g_j))$, $j=1,2$. Using Lemma \ref{ah}
we can find an integer $n_1 \in \N$ such that the element
$g_3=g_2g_1^{n_1} \in A$ is hyperbolic relatively to
$\mathfrak{C}_1$ and has infinite order. Applying Lemma \ref{Eg}
again, we achieve hyperbolicity of $G$ relative to
$\mathfrak{C}_2=\mathfrak{C}_1 \cup \{E_G(g_3)\}$. Set
$\mathcal{H}'= \bigsqcup_{H \in \mathfrak{C}_2} (H \setminus
\{1\})$.

 Let $\Omega \subset G$ be the finite subset and $K>0$ the constant
chosen according to Lemma~\ref{Omega} (where $G$ is considered to be relatively hyperbolic with respect to
$\mathfrak{C}_2$). Using Lemma \ref{ah} two more times, we can find  numbers
$m_1,m_2,m_3 \in \N$ such that
\begin{equation} \label{eq:choose_m}
g_i^{m_i} \notin \{y \in \langle \Omega \rangle \; : \; |y|_{\Omega} \le 21K\},\;\;\; i=1,2,3,
\end{equation}
and $h=g_1^{m_1}g_3^{m_3}g_2^{m_2}\in A$ is a hyperbolic element
(with respect to $\mathfrak{C}_2$) and has infinite order. Indeed,
first we choose $m_1$ to satisfy \eqref{eq:choose_m}. By Lemma
\ref{ah}, there is $m_3$ satisfying \eqref{eq:choose_m}, so that
$g_1^{m_1}g_3^{m_3} \in A^0$. Similarly $m_2$ can be chosen
sufficiently big to satisfy \eqref{eq:choose_m} and
$g_1^{m_1}g_3^{m_3}g_2^{m_2}\in A^0$. In particular, $h$ will be
non-commensurable with $g_j$, $j=1,2$ (otherwise, there would exist
$f \in G$ and $n \in \N$ such that $f^{-1}h^nf \in E(g_j)$, implying
$h \in fE(g_j)f^{-1}$ by Lemma \ref{Eg} and contradicting the hyperbolicity of $h$).

Consider a path $q$ labelled by the word $(g_1^{m_1} g_3^{m_3}
g_2^{m_2})^l$ in {\Ghh} for some $l \in \Z \setminus \{0\}$, where
each $g_i^{m_i}$ is treated as a single letter from
$\mathcal{H}'$. After replacing $q$ with $q^{-1}$, if necessary,
we assume that $l\in \N$. Let $p_1,\dots,p_{3l}$ be all components of $q$; by the construction of $q$, we
have $l(p_j)=1$ for each $j$. Suppose not all of these components
are isolated. Then one can find indices $1 \le s <t \le 3l$ and $i
\in \{1,2,3\}$ such that $p_s$ and $p_t$ are $E_G(g_i)$-components
of $q$, $(p_t)_-$ and $(p_s)_+$ are connected by a path $r$ with
$\phi(r) \in E_G(g_i)$, $l(r)\le 1$, and $(t-s)$ is minimal with
this property. To simplify the notation, assume that $i=1$ (the
other two cases are similar). Then  $p_{s+1}, p_{s+4}, \dots,
p_{t-2}$ are isolated $E_G(g_3)$-components of the cycle
$p_{s+1}p_{s+2}\dots p_{t-1} r$, and there are exactly $(t-s)/3
\ge 1$ of them. Applying Lemma \ref{Omega}, we obtain $g_3^{m_3}
\in \langle \Omega \rangle$ and $$\frac{t-s}{3} |g_3^{m_3}|_\Omega
\le K(t-s).$$ Hence $|g_3^{m_3}|_\Omega \le 3K,$ contradicting
\eqref{eq:choose_m}. Therefore two distinct components of $q$ can
not be connected with each other; that is, the path $q$ is
without backtracking.

To finish the proof of Lemma
\ref{lem:E^+_G} we need an auxiliary statement below.
Denote by $\mathcal{W}$ the set of all subwords of words
$(g_1^{m_1} g_3^{m_3} g_2^{m_2})^l$, $l \in \Z$ (where $g_i^{\pm
m_i}$ is treated as a single letter from $\mathcal H^\prime $).
Consider an arbitrary cycle $o=rqr'q'$ in {\Ghh}, where $\phi(q),
\phi(q') \in \mathcal{W}$; and set $C=\max\{ l(r),l(r')\}$. Let
$p$ be a component of $q$ (or $q'$). We will say that $p$ is {\it
regular} if it is not an isolated component of $o$. As $q$ and
$q'$ are without backtracking, this means that $p$ is either
connected to some component of $q'$ (respectively $q$), or to a
component of $r$, or $r'$.

\begin{lem}\label{lem:new_way} In the above notations
\begin{itemize}
\item[\rm (a)] if $C\le 1$ then every component of $q$ or $q'$ is regular;
\item[\rm (b)] if $C\ge 2$ then each of $q$ and $q'$ can have at most $15C$ components which are not regular.
\end{itemize}
\end{lem}

\begin{proof} Assume the contrary to (a). Then one can choose a cycle $o=rqr'q'$  with
$l(r),l(r') \le 1$, having at least one $E(g_i)$-isolated
component on $q$ or $q'$ for some $i\in \{1,2,3\}$, and such that
$l(q)+l(q')$ is minimal. Clearly the latter condition implies that
each component of $q$ or $q'$ is an isolated component of $o$.
Therefore $q$ and $q'$ together contain $k$ distinct $E(g_i)$-components of
$o$ where $k \ge 1$ and $k\ge \lfloor l(q)/3 \rfloor+\lfloor
l(q')/3 \rfloor$. Applying Lemma \ref{Omega} we obtain $g_i^{m_i}
\in \langle \Omega \rangle$ and $k|g_i^{m_i}|_\Omega \le
K(l(q)+l(q')+2)$, therefore $|g_i^{m_i}|_\Omega \le 11 K$,
contradicting the choice of $m_i$ in \eqref{eq:choose_m}.

Let us prove (b). Suppose that $C \ge 2$ and $q$ contains more than $15C$ isolated components of $o$.
 We consider two cases:

{\bf Case 1}. No component of $q$ is connected to a component of
$q'$. Then a component of $q$ or $q'$ can be regular only if it is
connected to a component of $r$ or $r'$. Since $q$ and $q'$ are
without backtracking, two distinct components of $q$ or $q'$ can
not be connected to the same component of $r$ (or $r'$). Hence $q$
and $q'$ together can contain at most $2C$ regular components.
Thus there is an index $i \in \{1,2,3\}$ such that the cycle $o$
has $k$ isolated $E(g_i)$-components, where $k \ge \lfloor l(q)/3
\rfloor + \lfloor l(q')/3 \rfloor -2C \ge \lfloor 5C \rfloor-2C
>2C>3$. By Lemma \ref{Omega}, $g_i^{m_i} \in \langle \Omega
\rangle$ and $k|g_i^{m_i}|_\Omega \le K(l(q)+l(q')+2C)$, hence
$$|g_i^{m_i}|_\Omega \le K\frac{3(\lfloor l(q)/3 \rfloor +1) +
3(\lfloor l(q')/3\rfloor+1)+2C}{\lfloor l(q)/3 \rfloor + \lfloor
l(q')/3 \rfloor -2C} \le K \left(3+\frac{6+8C}{2C} \right) \le 9
K,$$ contradicting the choice of $m_i$ in
\eqref{eq:choose_m}.

{\bf Case 2.} The path $q$ has at least one component which is
connected to a component of $q'$. Let $p_1,\dots,p_{l(q)}$ denote
the sequence of all components of $q$. By part (a), if $p_{s}$ and
$p_{t}$, $1 \le s \le t \le l(q)$, are connected to components of
$q'$, then for any $j$, $s \le j \le t$, $p_j$ is regular. We can
take $s$ (respectively $t$) to be minimal (respectively maximal)
possible. Consequently $p_1,\dots,p_{s-1}, p_{t+1},\dots,p_{l(q)}$
will contain the set of all isolated components of $o$ that belong
to $q$.

Without loss of generality we may assume that $s-1 \ge 15C/2$.
Since $p_s$ is connected to some component $p'$ of $q'$, there
exists a path $v$ in {\Ghh} satisfying $v_-=(p_{s})_-$,
$v_+=p'_+$, $\phi(v) \in \mathcal{H}'$, $l(v)=1$. Let $\bar q$
(respectively $\bar q'$) denote the subpath of $q$ (respectively
$q'$) from $q_-$ to $(p_s)_-$ (respectively from $p'_+$ to
$q'_+$). Consider a new cycle $\bar o = r \bar q v \bar q'$.
Reasoning as before, we can find $i \in \{1,2,3\}$ such that $\bar
o$ has $k$ isolated $E(g_i)$-components, where $k \ge \lfloor
l(\bar q)/3 \rfloor + \lfloor l({\bar q}')/3 \rfloor -C-1 \ge
\lfloor 15C/6 \rfloor-C-1 >C-1\ge 1$. Using Lemma \ref{Omega}, we
get $g_i^{m_i} \in \langle \Omega \rangle$ and
$k|g_i^{m_i}|_\Omega \le K(l(\bar q)+l(\bar q')+C+1)$. The latter
inequality implies $|g_i^{m_i}|_\Omega \le 21K$, yielding a
contradiction in the usual way and proving (b) for $q$.
By symmetry this property holds for $q'$ as well.
\end{proof}

Continuing the proof of Lemma \ref{lem:E^+_G}, consider an element
$x \in E_G(h)$. According to Lemma \ref{Eg}, there exists $l \in
\N$ such that
\begin{equation} \label{eq:xhx} xh^lx^{-1}=h^{\epsilon l},
\end{equation} where $\epsilon=\pm 1$.
Set $C=|x|_{X\cup \mathcal{H}'}$. After raising both sides of
\eqref{eq:xhx} in an integer power, we can assume that $l$ is
sufficiently large to satisfy $l>32C+3$.

Consider a cycle $o=rqr'q'$ in {\Ghh} satisfying $r_-=q'_+=1$, $r_+=q_-=x$,
$q_+=r'_-=xh^l$, $r'_+=q'_-=xh^lx^{-1}$, $\phi(q) \equiv (g_1^{m_1} g_3^{m_3} g_2^{m_2})^l$,
$\phi(q') \equiv (g_1^{m_1} g_3^{m_3} g_2^{m_2})^{-\epsilon l}$, $l(q)=l(q')=3l$, $l(r)=l(r')=C$.

Let $p_1,p_2,\dots,p_{3l}$ and $p'_1,p'_2,\dots,p'_{3l}$ be all components of $q$ and $q'$ respectively.
Thus, $p_3,p_6,p_9,\dots,p_{3l}$ are all $E_G(g_2)$-components of $q$. Since $l>17C$ and $q$ is without backtracking,
by Lemma \ref{lem:new_way}, there exist indices $1 \le s,s' \le 3l$ such that the $E_G(g_2)$-component $p_s$
of $q$ is connected to the $E_G(g_2)$-component $p'_{s'}$  of $q'$.
Without loss of generality, assume that $s \le 3l/2$ (the other
situation is symmetric). There is a path $u$ in {\Ghh} with $u_-=(p'_{s'})_-$, $u_+=(p_s)_+$, $\phi(u) \in E_G(g_2)$
and $l(u)\le 1$. We obtain a new cycle $o' =up_{s+1} \dots p_{3l}r'p'_1\dots p'_{s'-1}$ in the Cayley graph {\Ghh}.
Due to the choice of $s$ and $l$, the same argument as before will demonstrate that there are
$E_G(g_2)$-components $p_{\bar s}$, $p'_{\bar s'}$ of $q$, $q'$ respectively, which are connected and
$s<\bar s \le 3l$, $1 \le \bar s'<s'$ (in the case when $s>3l/2$, the same inequalities can be achieved by simply renaming
the indices correspondingly).

It is now clear that there exist $i \in \{1,2,3\}$  and connected $E_G(g_i)$-components
$p_{t}$, $p'_{t'}$ of $q$, $q'$ ($s <t\le 3l$, $1\le t' <s'$)  such that $t>s$ is minimal.
Let $v$ denote a path in {\Ghh} with $v_-=(p_{t})_-$, $v_+=(p_{t'})_+$, $\phi(v) \in E_G(g_i)$
and $l(v) \le 1$. Consider a  cycle $o''$ in  {\Ghh} defined by $o''=up_{s+1} \dots p_{t-1}vp'_{t'+1}\dots p'_{s'-1}$.
By part a) of Lemma \ref{lem:new_way}, $p_{s+1}$ is a regular component of the path $p_{s+1} \dots p_{t-1}$ in
$o''$ (provided that $t-1\ge s+1$). Note that $p_{s+1}$ can not be connected to $u$ or $v$ because $q$
is without backtracking, hence it must be
connected to a component of the path $p'_{t'+1}\dots p'_{s'-1}$. By the choice of $t$, we have $t=s+1$ and $i=1$.
Similarly $t'=s'-1$. Thus $p_{s+1}=p_t$ and
$p'_{s'-1}=p'_{t'}$ are connected $E_G(g_1)$-components of $q$ and $q'$.

In particular, we have $\epsilon=1$. Indeed, otherwise we would have
$\phi(p_{s'-1}) \equiv g_3^{m_3}$ but $g_3^{m_3} \notin E_G(g_1)$.
Therefore $x \in E^+_G(h)$ for any $x \in E_G(h)$, consequently
$E_G(h)=E^+_G(h)$.

Observe that $u_-=v_+$ and $u_+=v_-$, hence $\phi(u)$ and $\phi(v)^{-1}$ represent the same element
$z \in E_G(g_2) \cap E_G(g_1)$. By construction, $x=h^{\alpha}zh^{\beta}$ where $\alpha=(3l-s')/3 \in \Z$, and
$\beta =-s/3 \in \Z$. Thus
$x \in \langle h, E_G(g_1) \cap E_G(g_2) \rangle$ and the first part of the claim 2 is proved.

Assume now that $E_G(g_j)=E^+_G(g_j)$ for $j=1,2$. Then
$h=g_1^{m_1}(g_2g_1^{n_1})^{m_3}g_2^{m_2}$ belongs to the centralizer of the
finite subgroup $E_G(g_1) \cap E_G(g_2)$ (because of the choice of
$g_1,g_2$ above). Consequently $E_G(h)=\langle h \rangle \times
(E_G(g_1) \cap E_G(g_2))$.
\end{proof}

\begin{lem} \label{lem:good-elem} Let $S$ be a non-elementary
subgroup of $G$ with $S^0 \neq \emptyset$. Then
\begin{enumerate}
\item[\rm(i)] there exist non-commensurable elements $h_1, h'_1 \in S^0$
with $E_G(h_1) \cap E_G(h_1')=E_G(S)$;
\item[\rm(ii)] $S^0$ contains an element $h$ such that $E_G(h)=\langle h
\rangle \times E_G(S)$.
\end{enumerate}
\end{lem}

\begin{proof} Choose an element $g_1 \in S^0$. By Lemma \ref{Eg}, $G$ is hyperbolic relative to the
collection $\mathfrak{C}=\Hl \cup \{E_G(g_1)\}$. Since the
subgroup $S$ is non-elementary, there is $a \in S \setminus
E_G(g_1)$, and Lemma \ref{ah} provides us with an integer  $n \in
\N$ such that $g_2=ag_1^{n} \in S$ is a hyperbolic element of
infinite order (now, with respect to the family of peripheral
subgroups $\mathfrak{C}$). In particular, $g_1$ and $g_2$ are
non-commensurable and hyperbolic relative to $\Hl$.

Applying Lemma \ref{lem:E^+_G}, we   find $h_1 \in S^0$ (with
respect to the collection of peripheral subgroups $\Hl$) with
$E_G(h_1)=E^+_G(h_1)$ such that $h_1$ is not commensurable with
$g_j$, $j=1,2$. Hence, $g_1$ and $g_2$ stay hyperbolic after
including $E_G(h_1)$ into the family of peripheral subgroups (see
Lemma \ref{Eg}). This allows to construct (in the same manner) one
more element $h_2 \in \langle g_1,g_2\rangle \le S$ which is
hyperbolic relative to $(\Hl \cup E_G(h_1))$ and satisfies
$E_G(h_2)=E^+_G(h_2)$. In particular, $h_2$ is not commensurable
with $h_1$.

We claim now that there exists $x \in S$ such that
$E_G(x^{-1}h_2x) \cap E_G(h_1)=E_G(S)$.  By definition, $E_G(S)
\subseteq E_G(x^{-1}h_2x) \cap E_G(h_1)$.  To obtain the inverse
inclusion, arguing by the contrary, suppose that for each $x \in
S$ we have
\begin{equation} \label{eq:inter-e} (E_G(x^{-1}h_2x) \cap E_G(h_1)) \setminus E_G(S) \neq \emptyset.\end{equation}
Note that  if $g\in S^0$ with $E_G(g)=E^+_G(g)$, then the set of
all elements of finite order in $E_G(g)$ form a finite subgroup
$T(g) \le E_G(g)$ (this is a well-known property of groups, all of
whose conjugacy classes are finite). The elements $h_1$ and $h_2$
are not commensurable, therefore
$$E_G(x^{-1}h_2x) \cap E_G(h_1) =T(x^{-1}h_2x) \cap T(h_1)=x^{-1}T(h_2)x \cap T(h_1).$$
For each pair of elements $(b,a) \in D=T(h_2) \times
(T(h_1)\setminus E_G(S))$ choose $x=x(b,a) \in S$ so that
$x^{-1}bx=a$ if such $x$ exists; otherwise set $x(b,a)=1$.

The assumption \eqref{eq:inter-e} clearly implies that
$\displaystyle S=\bigcup_{(b,a) \in D} x(b,a) C_S(a)$, where
$C_S(a)$ denotes the centralizer of $a$ in $S$. Since the set $D$
is finite, a well-know theorem of B. Neumann \cite{Neumann}
implies that there exists $a \in T(h_1) \setminus E_G(S)$ such
that $|S:C_S(a)|<\infty$. Consequently, $a \in E_G(g)$ for every
$g \in S^0$,  that is, $a \in E_G(S)$, a contradiction.

Thus, $E_G(xh_2x^{-1}) \cap E_G(h_1)=E_G(S)$ for some $x \in S$.
After setting  $h'_1=x^{-1}h_2x \in S^0$, we see that elements
$h_1$ and  $h'_1$ satisfy the claim (i). Since
$E_G(h'_1)=x^{-1}E_G(h_2)x$, we have $E_G(h'_1)=E^+_G(h'_1)$. To
demonstrate (ii), it remains to apply Lemma \ref{lem:E^+_G} and
obtain an element  $h \in \langle h_1,h'_1 \rangle \le S$ which
has the desired properties.
\end{proof}

\begin{proof}[Proof of Proposition \ref{suit}]
The implication  $(1) \Rightarrow (2)$ is an immediate consequence
of the definition. The inverse implication follows directly from
the first claim of Lemma \ref{lem:good-elem} ($S$ is
non-elementary as $S^0 \neq \emptyset$ and
$E_G(S) =\{1\}$).
\end{proof}

\section{Proofs of the main results}

The following simplification of Theorem 2.4 from \cite{SCT} is the
key ingredient of the proofs in the rest of the paper.

\begin{thm}\label{glue}
Let $U$ be a group hyperbolic relative to a collection of subgroups
$\Vl $, $S$ a suitable subgroup of $U$, and $T$ a finite subset of
$U$. Then there exists an epimorphism $\eta \colon U\to W$ such
that:
\begin{enumerate}
\item The restriction
of $\eta $ to  $\bigcup_{\lambda \in \Lambda} V_\lambda $ is injective, 
and the group $W$ is properly relatively hyperbolic with respect to the
collection $\{\eta (V_\lambda )\}_{\lambda \in \Lambda }$.

\item For every $t\in T$, we have $\eta (t)\in \eta (S)$.

\end{enumerate}
\end{thm}

Let us also mention two  known results we will use. The first
lemma is a particular case of Theorem 1.4 from  \cite{RHG}  (if $g
\in G$ and $H \le G$, $H^g$ denotes the conjugate $g^{-1}Hg \le
G$).

\begin{lem}\label{malnorm}
Suppose that a group $G$ is hyperbolic relative to a collection of
subgroups $\Hl $. Then
\begin{enumerate}
\item[(a)]  For any $g\in G$ and any $\lambda , \mu \in \Lambda $,
$\lambda \ne \mu $, the intersection $H_\lambda^g\cap H_\mu $ is
finite.
\item[(b)] For any $\lambda \in \Lambda $ and any $g\notin
H_\lambda $, the intersection $H_\lambda^g \cap H_\lambda $ is
finite.
\end{enumerate}
\end{lem}

The second result can easily be derived from Lemma \ref{ah}.

\begin{lem}[Corollary 4.5, \cite{ESBG}] \label{prop}
Let $G$ be an infinite properly relatively hyperbolic group. Then
$G$ contains a hyperbolic element of infinite order.
\end{lem}

\begin{lem}\label{lem:qbf} Let the group $G$ be hyperbolic with respect to the collection
of peripheral subgroups $\Hl$ and let $N \lhd G$ be a finite normal subgroup. Then

\begin{enumerate}
\item If $H_\lambda $ is infinite for some $\lambda \in \Lambda $,
then $N\le H_\lambda $;
\item The quotient $\bar G=G/N$ is hyperbolic relative to the natural
image of the collection $\Hl $.
\end{enumerate}
\end{lem}

\begin{proof}
Let $K_\lambda $, $\lambda \in \Lambda $, be the kernel of the
action of $H_\lambda $ on $N$ by conjugation. Since $N$ is finite,
$K_\lambda $ has finite index in $H_\lambda $. On the other hand
$K_\lambda \le H_\lambda \cap H_\lambda ^g$ for every $g\in N$. If
$H_\lambda $ is infinite this implies $N\le H_\lambda $ by Lemma
\ref{malnorm}.

To prove the second assertion, suppose that $G$ has a relatively
finite presentation \eqref{G} with respect to the free product $F$
defined in \eqref{F}. Denote by $\bar X$ and $\bar H_\lambda$ the
natural images of $X$ and $H_\lambda$ in $\bar G$.
In order to show that $\bar G$ is relatively hyperbolic, one has to consider it as a quotient 
of the free product $\bar F=(\ast_{\lambda \in \Lambda} \bar
H_\lambda)*F(\bar X)$. As $G$ is a quotient of $F$, we can choose
some finite preimage $M\subset F$ of $N$. For each element $f \in
M$, fix a word in $X \cup \mathcal{H}$ which represents it in $F$
and denote by $\mathcal{S}$ the (finite) set of all such words. By
the universality of free products, there is a natural epimorphism
$\varphi: F \to \bar F$ mapping $X$ onto $\bar X$ and each
$H_\lambda$ onto $\bar H_\lambda$. Define the subsets
$\bar{\mathcal{R}}$ and $\bar{\mathcal{S}}$ of words in $\bar X
\cup \bar{\mathcal{H}}$ (where $\bar{\mathcal{H}} =
\bigsqcup_{\lambda \in \Lambda} (\bar H_\lambda \setminus \{1\})$)
by $\bar{\mathcal{R}}=\varphi(\mathcal R)$ and
$\bar{\mathcal{S}}=\varphi(\mathcal{S})$. Then the group $\bar G$
possesses the  relatively  finite presentation
\begin{equation} \label{eq:G_1} \langle \bar X,\; \{ \bar H_\lambda\}_{\lambda\in \Lambda} \mid \bar R=1,\,
\bar R\in\bar{\mathcal{R}};\,  \bar S=1,\,\bar S\in\bar{\mathcal{S}}\rangle.\end{equation}

Let $\psi: F \to G$ denote the natural epimorphism and
$D=\max\{\|s\|~:~s \in \mathcal{S}\}$. Consider any non-empty word
$\bar w$ in the alphabet $\bar X \cup \bar{\mathcal{H}}$
representing the identity in $\bar G$. Evidently we can choose a
word $w$ in $X \cup \mathcal{H}$ such that $\bar w =_{\bar F}
\varphi(w)$ and $\|w\|=\|\bar w\|$. Since $\ker(\psi) \cdot M$ is
the kernel of the induced homomorphism from $F$ to $\bar G$, we
have $w=_{F}vu$ where $u \in \mathcal{S}$ and $v$ is a word in $X
\cup \mathcal{H}$ satisfying $v=_{G} 1$ and $\|v\| \le \|w\|+D$.
Since $G$ is relatively hyperbolic there is a constant $C \ge 0$
(independent of $v$) such that
$$ v=_F\prod\limits_{i=1}^k f_i^{-1}R_i^{\pm 1}f_i,$$ where $R_i \in \mathcal{R}$, $f_i \in
F$, and $k \le C \|v\| $.  Set $\bar R_i=\varphi(R) \in
\bar{\mathcal{R}}$, $\bar f_i=\varphi(f_i) \in \bar F$,
$i=1,2,\dots,k$, and ${\bar R}_{k+1}=\varphi(u) \in
\bar{\mathcal{S}}$, ${\bar f}_{k+1} =1$. Then $$\bar w =_{\bar F}
\prod \limits_{i=1}^{k+1} {\bar f}_i^{-1} \bar R_i^{\pm 1} {\bar
f}_i, $$ where
$$k+1 \le C\|v\|+1 \le C(\|w\|+D)+1 \le C\|\bar w\| +CD + 1 \le (C+CD+1)\|\bar w\|.$$
Thus, the relative presentation \eqref{eq:G_1} satisfies a linear
isoperimetric inequality with the constant $(C+CD+1)$.
\end{proof}

Now we are ready to prove Theorem \ref{SQ}.

\begin{proof}[Proof of Theorem \ref{SQ}]
Observe that the quotient of $G$ by the finite normal subgroup
$N=E_G(G)$ is obviously non-elementary.
Hence the image of any finite $H_\lambda $ is a proper subgroup of $G/N$.
On the other hand, if $H_\lambda$ is infinite, then $N \le H_\lambda \lneqq G$
by Lemma \ref{lem:qbf}, hence its image is also proper in $G/N$.
Therefore $G/N$ is properly relatively hyperbolic with respect to
the collection of images of $H_\lambda $, $\lambda \in \Lambda $
(see Lemma \ref{lem:qbf}).  Lemma \ref{Ashot's Lemma} implies
$E_{G/N}(G/N)=\{ 1\} $. Thus, without loss of generality, we may
assume that $E_G(G)=1$.

It is straightforward to see that the free product $U=G\ast R$ is
hyperbolic relative to the collection $\Hl \cup \{ R\} $ and
$E_{G\ast R} (G)=E_G(G)=1$. Note  that $G^0$ is non-empty by
Lemma \ref{prop}. Hence $G$ is a suitable subgroup of $G\ast R$ by
Proposition \ref{suit}. Let $Y$ be a finite generating set of $R$.
It remains to apply Theorem \ref{glue} to $U=G\ast R$, the obvious
collection of peripheral subgroups, and the finite set $Y$.
\end{proof}

 To prove Theorem \ref{comquot} we need
one more auxiliary result which was proved in the full generality in \cite{RHG} (see
also \cite{F}):

\begin{lem}[Theorem 2.40, \cite{RHG}]\label{exhyp}
Suppose that a group $G$ is hyperbolic relative to a collection of
subgroups $\Hl \cup \{ S_1, \ldots , S_m\} $, where $S_1, \ldots, S_m $
are hyperbolic in the ordinary (non-relative) sense. Then $G$ is hyperbolic relative to $\Hl $.
\end{lem}

\begin{proof}[Proof of Theorem \ref{comquot}]
 Let $G_1$, $G_2$ be finitely generated groups which are properly relatively hyperbolic with respect to
collections of subgroups $\{ H_{1\lambda }\} _{ \lambda \in
\Lambda } $ and $\{ H_{2\mu }\} _{ \mu \in M }$ respectively.
Denote by $X_i$ a finite generating set of the group $G_i$, $i=1,2
$. As above we may assume that $E_{G_1}(G_1)=E_{G_2}(G_2)=\{ 1\}
$. We set $G=G_1\ast G_2$. Observe that $E_G(G_i)=E_{G_i}(G_i)=\{
1\} $ and hence $G_i$ is suitable in $G$ for $i=1,2$ (by Lemma \ref{prop} and
Proposition \ref{suit}).

By the definition of suitable subgroups, there are two
non-commensurable elements  $g_1, g_2\in G_2^0$ such that
$E_G(g_1)\cap E_G(g_2)=\{ 1\} $. Further, by Lemma \ref{Eg}, the
group $G$ is hyperbolic relative to the collection $\mathfrak P=
\{ H_{1\lambda }\} _{ \lambda \in \Lambda } \cup \{ H_{2\mu }\} _{
\mu \in M } \cup \{ E_G(g_1), E_G(g_2)\} $. We now apply Theorem
\ref{glue} to the group $G$ with the collection of peripheral
subgroups $\mathfrak P$, the suitable subgroup $G_1\le G$, and the
subset $T=X_2$.  The resulting group $W$ is obviously a quotient
of $G_1$.

Observe that $W$ is hyperbolic relative to (the image of) the
collection $\{ H_{1\lambda }\} _{ \lambda \in \Lambda } \cup \{
H_{2\mu }\} _{ \mu \in M } $ by Lemma \ref{exhyp}. We  would like
to show that $G_2$ is a suitable subgroup of $W$ with respect to
this collection. To this end we note that $\eta (g_1)$ and $\eta
(g_2) $ are elements of infinite order as $\eta $ is injective on
$E_G(g_1)$ and $E_G(g_2)$.  Moreover, $\eta (g_1)$ and $\eta (g_2)
$ are not commensurable in $W$.  Indeed, otherwise, the
intersection $\bigr(\eta(E_G(g_1))\bigl)^g \cap \eta (E_G(g_2))$
is infinite for some $g\in G$ that contradicts the first assertion
of Lemma \ref{malnorm}. Assume now that $g\in E_{W}(\eta (g_i))$
for some $i\in \{ 1,2\} $. By the first assertion of Lemma
\ref{Eg}, $\big( \eta (g_i^m)\big) ^g= \eta (g_i^{\pm m})$ for
some $m\ne 0$. Therefore, $\big( \eta (E_G(g_i)) \big) ^g\cap \eta
(E_G(g_i))$ contains $\eta (g_i^m)$  and, in particular, this
intersection is infinite. By the second assertion of Lemma
\ref{malnorm}, this means that $g\in \eta (E_G(g_i))$.  Thus,
$E_{W}(\eta (g_i))=\eta (E_G(g_i))$. Finally, using injectivity of
$\eta $ on  $E_G(g_1) \cup E_G(g_2)$, we obtain
$$ E_{W}(\eta (g_1))\cap E_{W}(\eta (g_2))=\eta (E_G(g_1))\cap
\eta (E_G(g_2))=\eta \big( E_G(g_1)\cap E_G(g_2)\big) = \{ 1\} .$$
This means that the image of $G_2$ is a suitable subgroup of $W$.

Thus we may apply Theorem \ref{glue} again to the group $W$, the
subgroup $G_2$ and the finite subset $X_1$. The resulting group $Q$
is the desired common quotient of $G_1$ and $G_2$. The last property, which claims
that $Q$ can be obtained from $G_1 \ast G_2$ by adding only finitely many relations,
follows because $G_1\ast G_2$ and $G$ are hyperbolic with respect to the same family of
peripheral subgroups and any relatively hyperbolic group is relatively finitely presented.
\end{proof}

\small \textsc{
G. Arzhantseva, Universit\'{e} de Gen\`{e}ve,
Section de Math\'{e}matiques,
2-4 rue du Li\`{e}vre,
Case postale 64,
1211 Gen\`{e}ve 4, Switzerland}

{\it Email:} \texttt{Goulnara.Arjantseva@math.unige.ch}

\vspace{.3cm}
\textsc{
A. Minasyan, Universit\'{e} de Gen\`{e}ve,
Section de Math\'{e}matiques,
2-4 rue du Li\`{e}vre,
Case postale 64,
1211 Gen\`{e}ve 4, Switzerland}

{\it Email:} \texttt{aminasyan@gmail.com}

\vspace{.3cm} \textsc{D. Osin, NAC 8133, Department of
Mathematics, The City College of the City University of New York,
Convent Ave. at 138th Street, New York, NY 10031, USA}

{\it Email:} \texttt{denis.osin@gmail.com}
\end{document}